\documentclass[12pt,a4paper,oneside,titlepage]{article}
\usepackage{amsmath,amsfonts,amssymb,amsthm,amscd}
\usepackage{graphicx,xcolor}
\newtheorem{Theorem}{Theorem}[section]

\newcommand{\R}{\mathbb R}

\newcommand{\Z}{\mathbb Z}
\newcommand{\N}{\mathbb N}

\title{On a Theorem of Nathanson on Diophantine Approximation}
\author{Jaroslav Han\v cl
and Tho Phuoc Nguyen}

\begin{document}

\maketitle
\date{}

\begin{abstract}

In 1974, M. B. Nathanson proved that every irrational number $\alpha$ represented by a simple continued fraction with infinitely many elements greater than or equal to $k$ is approximable by an infinite number of rational numbers $p/q$ satisfying $|\alpha-p/q|<1/(\sqrt{k^2+4}q^2)$. In this paper we refine this result.
\end{abstract}

\section{Introduction}

For $k\geq1$, let $F(k)$ denote the set of all real numbers $x$ such that $0\leq x \leq1$ and the simple continued fraction for $x$ has no partial quotient greater than $k$.
Following Dirichlet \cite{dirichlet}, Hurwitz \cite{hurwitz}, Borel \cite{borel1}, \cite{borel2}, Vahlen \cite{vahlen} and Cohn \cite{cohn},   in 1974, Nathanson \cite{nathanson}, proved the following. 

\begin{Theorem} Suppose $k\geq1$ and let $x$ be a real irrational number not equivalent to an element of $F(k-1)$. Then there are infinitely many rational numbers $\dfrac{p}{q}$ such that
$$
\left|\alpha-\dfrac{p}{q} \right|<\dfrac{1}{q^2\sqrt{k^2+4}}.
$$
The constant $\dfrac{1}{\sqrt{k^2+4}}$ is best possible.
\end{Theorem}
In this paper we refine this theorem, replacing $\sqrt{k^2+4}$ by the function 
 $$f(q)=\dfrac{q^2\sqrt{k^2+4}}{2}\left(1+\sqrt{1+\dfrac{4}{(k^2+4)q^2}} \right).$$   
These results have a history. 
Dirichlet \cite{dirichlet} showed that if $x$ is an irrational number then there exist infinitely rational numbers $\frac pq$ such that 
$\mid x-\frac pq\mid <\frac 1{q^2}=g(q)$. 
In 1891 Hurwitz \cite{hurwitz} improved this result to  $g(q)=\frac 1{\sqrt{5}q^2}$  and in 2015 Han\v cl \cite{hancl1} further improved this to
$g(q)=\left(q^2\frac{\sqrt 5}2\left (1+\sqrt{1+\frac 4{5q^2}}\right)\right)^{-1}$. In 1973, Cohn \cite{cohn} had also presented a simple proof of the Hurwitz's theorem.

In 1903, 
Borel 
\cite{borel1}, \cite{borel2} proved that 
if $\alpha\in\R$, $n\in\N$ and $\frac{p_{n-1}}{q_{n-1}}$, $\frac{p_n}{q_n}$ and $\frac{p_{n+1}}{q_{n+1}}$ are three consecutive convergents of the number $\alpha$, 
then at least one of them satisfies the inequality
$
\Bigl\lvert \alpha -\frac pq\Bigr\rvert< \frac 1 {\sqrt 5 q^2}\,.
$
The constant  $\sqrt 5$ can be replaced  by  $\sqrt 5+\frac {4-5\sqrt 5+\sqrt{61}}{2q^2}$\cite{hanclnair}.

An important precursor result to Borel's theorem is the theorem of Vahlen \cite{vahlen}\label{hb1t.1}  from 1895. He proved, for a real number $\alpha$, a positive integer $n$ and two consecutive convergents $\frac{p_{n-1}}{q_{n-1}}$, $\frac{p_n}{q_n}$  of the number $\alpha$,
then at least one of them satisfies the inequality
$
\Bigl\lvert \alpha -\frac pq\Bigr\rvert< \frac 1 {2q^2}\,.
$
Other related results concerning  Vahlen's  theorem can be found in  \cite{bh}, \cite{hancl1} or \cite{hancl2}.

An excellent source of basic background is  Hardy and Wright \cite{hardy}. The books of  Wall \cite{wall}, Hensley \cite{hensley}, Khinchin \cite{khinchin}, Karpenkov \cite{karpenkov}, Thron \cite{jones},  Rosen \cite{rosen}, Borwein and Borwein \cite{borwein} and Feldman and Nesterenko \cite{feldman} are also very useful.

We now formally state our main result improving Nathanson's result. 
\begin{Theorem} \label{nathanson.t1}  
Let $\alpha=[a_0;a_1,\dots]$ be an irrational number not equivalent to an element of $F(k-1)$. Then there are infinitely many integers $p$ and $q$ such that
\begin{equation} \label{eq.theorem2}
    \left|\alpha-\dfrac{p}{q}\right| \leq \dfrac{1}{\dfrac{q^2\sqrt{k^2+4}}{2}\left(1+\sqrt{1+\dfrac{4}{(k^2+4)q^2}} \right)}
\end{equation}

The function $f(q)=\dfrac{q^2\sqrt{k^2+4}}{2}\left(1+\sqrt{1+\dfrac{4}{(k^2+4)q^2}} \right)$ is best possible in the sense that there isn't any function $g(q)$ such that $g(q)>f(q)$ for all positive integers $q$ and $g(q)$ satisfies \eqref{eq.theorem2} for infinitely many integers $p$ and $q$.

The equality in \eqref{eq.theorem2} is obtained only if the number $\alpha=\alpha_1=\dfrac{\sqrt{k^2+4}-k}{2}=[0;\overline{k}]$ or $\alpha=\alpha_2=\dfrac{k+2-\sqrt{k^2+4}}{2}=[0;1,k-1,\overline{k}]$ plus an integer.
\end{Theorem}

\section{Notation} 

Throughout the paper, $\Z$, $\N$, $\N_0$ and $\R$ will denote the sets of integers, positive integers, non-negative integers and real numbers respectively. Let $\alpha$ be a real number and suppose $n\in\N_0$\,.  Let $\alpha=[a_0;a_1,a_2,\dots ]$ be  its simple continued fraction expansion.  Also let $\dfrac {p_n}{q_n}=[a_0;a_1,a_2,\dots ,a_n]$ be its $n$-th convergent. 
The following recurrence relations for convergents are known
\begin{align*}
p_0&=a_0\,, &  p_1&=a_1a_0+1\,, &  p_{n+2}&=a_{n+2}p_{n+1}+p_n\,,\\
q_0&=1\,, &    q_1&=a_1\,,   &          q_{n+2}&=a_{n+2}q_{n+1}+q_n\,.
\end{align*}
For a simple continued fraction expansion  we have that
$$a=[a_0;a_1,a_2,\dots ]=[a_0;a_1,a_2,\dots ,a_n,[a_{n+1};a_{n+2},a_{n+3},\dots ]].$$
Taking a difference of two consecutive convergents we obtain that
$$
q_{n+1}p_n-p_{n+1}q_n=(-1)^{n+1}\,.
$$ 
Finally we have the identity most frequently used in this article.
\begin{equation}\label{hbor2}
\Bigl\lvert a -\frac{p_n}{q_n}\Bigr\rvert =
  \frac 1{q_n^2([a_{n+1};a_{n+2},\dots]+[0;a_n,a_{n-1},\dots ,a_1])}\,,
\end{equation}
 where if $n=0$\,, then we set $[0;a_n,a_{n-1},\dots,a_1]=0$\,. 

For simple continued fraction expansions, if we have $a=[a_0;a_1,a_2,\dots ,a_k]$ for finite $k\geq 1$\,, then we suppose that $a_k\not= 1$\,.  More details on the discussion in this section can be found in \cite{schmidt}, pages $7$ to $10$\,.

If $\alpha=[0;\overline{k}]=\dfrac{\sqrt{k^2+4}-k}{2}$, then we have $p_0=0,p_1=q_0=1,q_1=k$, $p_{n+2}=kp_{n+1}+p_n$ and $q_{n+2}=kq_{n+1}+q_n$ for all $n\in \N_0$. 

If $\alpha=[0;1,k-1,\overline{k}]=\dfrac{k+2-\sqrt{k^2+4}}{2}$, then we have $p_0=0,q_0=q_1=p_1=1$, $p_2=k-1,q_2=k$, $p_{n+2}=kp_{n+1}+p_n$ and $q_{n+2}=kq_{n+1}+q_n$ for all $n \in \N$. For every $n \in \N$ we have
$$(-1)^n\left(\dfrac{\sqrt{k^2+4}-k}{2} - \dfrac{p_n}{q_n} \right)=\dfrac{1}{\dfrac{q_n^2\sqrt{k^2+4}}{2}\left( 1+ \sqrt{1+\dfrac{4(-1)^{n+1}}{(k^2+4)q_n^2}} \right)}$$
and
$$(-1)^n\left(\dfrac{k+2-\sqrt{k^2+4}}{2} - \dfrac{p_n}{q_n} \right)=\dfrac{1}{\dfrac{q_n^2\sqrt{k^2+4}}{2}\left( 1+ \sqrt{1+\dfrac{4(-1)^{n}}{(k^2+4)q_n^2}} \right)}$$
which follow from the fact that $$p_n=\dfrac{\sqrt{k^2+4}}{k^2+4}\left(\dfrac{2}{\sqrt{k^2+4}-k} \right)^n - \dfrac{\sqrt{k^2+4}}{k^2+4}\left(\dfrac{-2}{\sqrt{k^2+4}+k} \right)^n$$
and
\begin{align*}
    q_n=\dfrac{k^2+4+k\sqrt{k^2+4}}{2k^2+8}&\left(\dfrac{2}{\sqrt{k^2+4}-k} \right)^n \\
    & + \dfrac{k^2+4-k\sqrt{k^2+4}}{2k^2+8}\left(\dfrac{-2}{\sqrt{k^2+4}+k} \right)^n
\end{align*}
 for $\alpha=\dfrac{\sqrt{k^2+4}-k}{2}$. At the same time
 \begin{align*}
     p_n=\dfrac{k^2+4+(k-2)\sqrt{k^2+4}}{2k^2+8}&\left(\dfrac{2}{\sqrt{k^2+4}-k} \right)^{n-1} \\
     &+ \dfrac{k^2+4+(2-k)\sqrt{k^2+4}}{2k^2+8}\left(\dfrac{-2}{\sqrt{k^2+4}+k} \right)^{n-1}
 \end{align*}
 and 
 \begin{align*}
     q_n=\dfrac{k^2+4+k\sqrt{k^2+4}}{2k^2+8}&\left(\dfrac{2}{\sqrt{k^2+4}-k} \right)^{n-1} \\
     &+ \dfrac{k^2+4-k\sqrt{k^2+4}}{2k^2+8}\left(\dfrac{-2}{\sqrt{k^2+4}+k} \right)^{n-1}
 \end{align*}
 for $\alpha=\dfrac{k+2-\sqrt{k^2+4}}{2}$.
 
 All this can be proved by mathematical induction. From this, for\\ $\alpha=\dfrac{\sqrt{k^2+4}-k}{2}$ and $n$ odd or $\alpha=\dfrac{k+2-\sqrt{k^2+4}}{2}$ and $n$ even, we have
 \begin{equation} \label{worst}
     \left|\alpha - \dfrac{p_n}{q_n} \right|=\dfrac{1}{\dfrac{q_n^2\sqrt{k^2+4}}{2}\left( 1+ \sqrt{1+\dfrac{4}{(k^2+4)q_n^2}} \right)}.
 \end{equation}

\section{Proofs}

\begin{proof}[Proof of Theorem \ref{nathanson.t1}] 
We look for the worst approximation such that $a_n \geq k$ for infinitely many $n$. From \eqref{worst} and the fact that $\sqrt{1+x}<1+\dfrac{1}{2}x$ for all $x\in(0,1)$ we have
\begin{align*}
    \dfrac{\sqrt{k^2+4}}{2}\left( 1+\sqrt{1+\dfrac{4}{(k^2+4)q_n^2}} \right) &< \dfrac{\sqrt{k^2+4}}{2}\left( 1+1+\dfrac{1}{2}\dfrac{4}{(k^2+4)q_n^2} \right) \\
    &=\sqrt{k^2+4}+\dfrac{1}{\sqrt{k^2+4}q_n^2} \xrightarrow[]{n \to \infty} \sqrt{k^2+4}.
\end{align*}

Suppose that $\alpha$ is not equivalent to $\alpha_1$. We now consider some cases.
\begin{enumerate}
    \item 
    Assume that $a_{n+1}\geq k+2$ for infinitely many $n$. From this we obtain that
    \begin{align*}
        [a_{n+1};a_{n+2},\dots]+[0;a_n,\dots ,a_1] &> k+2>\sqrt{k^2+4}+\dfrac{1}{\sqrt{k^2+4}} \\
        &\geq\sqrt{k^2+4}+\dfrac{1}{\sqrt{k^2+4}q_n^2}.
    \end{align*}
    
    \item 
    Suppose that $a_n \in \{1,2,\dots,k+1\}$ for $n\geq m$. Hence there is infinitely many $n>m$ such that $a_{n+1}=k+1$. This implies that
\begin{align*}
    [a_{n+1}&;a_{n+2},\dots]+[0;a_n,\dots ,a_1] \\
    &=k+1+[0;a_{n+2},\dots]+[0;a_n,a_{n-1}\dots a_{m},a_{m-1}.\dots ,a_1]\\
    &\geq k+1+[0;\overline{k+1,1}]+[0;k+1,1,k+1,1,\dots,a_{m-1},\dots ,a_1].
\end{align*}
Let $H=k+1+[0;\overline{k+1,1}]+[0;k+1,1,k+1,1,\dots,a_{m-1},\dots ,a_1]$. Then we have that
$$\lim_{n \rightarrow \infty}{H}=k+1+2[0;\overline{k+1,1}]=\dfrac{k^2+k+\sqrt{k^2+6k+5}}{k+1}>\sqrt{k^2+4}.$$

    \item 
    Assume that $a_n \in \{1,2,\dots,k\}$ for $n\geq m$. Suppose that there are infinitely many $n\in \N$ and $i=i(n)\in \Z,i\geq -1$ such that \\
    $a_{n+1},a_n,a_{n-1},\dots,a_{n-i}=k$ and $a_{n+2},a_{n-i-1}<k$. We now consider two cases
    \begin{enumerate}
        \item Suppose that $i+2$ is odd. From this we obtain that
        \begin{align*}
            a_{n+1}&+[0;a_{n+2},\dots]+[0;a_n,a_{n-1},\dots,a_{n-i},a_{n-i-1},\dots,a_1] \\
            &=k+[0;a_{n+2},\dots]+[0;\underbrace{k,k,\dots,k}_{i+1\text{ elements}},a_{n-i-1},\dots,a_1]\\
            &\geq k+[0;k-1,a_{n+3}\dots]+[0;\underbrace{k,k,\dots,k}_{i+1\text{ elements}},k-1,\dots,a_1]\\
            &> k+[0;k-1,k+1]+[0;\underbrace{k,k,\dots,k}_{i+1\text{ elements}},k,\dots,a_1]\\
            &\xrightarrow[]{n \to \infty} \dfrac{k^3+2k+2+k^2\sqrt{k^2+4}}{2k^2}>\sqrt{k^2+4}.
        \end{align*}
         \item Now let $i+2$ be even. Set 
        $$A=[0;a_{n-i-1},a_{n-i-2},\dots,a_1]>\dfrac{1}{k}.$$
        $$B=[0;a_{n+2},a_{n+3},\dots]>\dfrac{1}{k}.$$
        \begin{enumerate}
            \item Suppose that $A>B$ and consider
            \begin{align*}
                &a_{n-i}+[0;a_{n-i+1},\dots,a_{n+2},\dots]+[0;a_{n-i-1},\dots,a_1]-\sqrt{k^2+4} \\
                &=k+[0;\underbrace{k,k,\dots,k}_{i+1\text{ elements}},a_{n+2},\dots]+[0;a_{n-i-1},\dots,a_1]-\sqrt{k^2+4} \\
                &\geq k+[0;k,a_{n+2},\dots]+[0;a_{n-i-1},\dots,a_1]-\sqrt{k^2+4} \\
                &= k+[0;k+B]+A-\sqrt{k^2+4} \\
                &>k+[0;k+B]+B-\sqrt{k^2+4}.
            \end{align*}
            From this and the fact that $f(B)=\dfrac{1}{k+B}+B$ is an increasing function for $B>\dfrac{1}{k}$, we obtain that for all $B$
            $$f(B)>f\left( \dfrac{1}{k}\right)=\dfrac{1}{k+\dfrac{1}{k}}+\dfrac{1}{k}.$$
            Hence
            \begin{align*}
                k+[0;k+B]+B&-\sqrt{k^2+4}>k+\dfrac{1}{k+\dfrac{1}{k}}+\dfrac{1}{k}-\sqrt{k^2+4}\\
                &=\dfrac{\left( k+\dfrac{1}{k+\dfrac{1}{k}}+\dfrac{1}{k} \right)^2-(k^2+4)}{k+\dfrac{1}{k+\dfrac{1}{k}}+\dfrac{1}{k}+\sqrt{k^2+4}}\\
                &=\dfrac{\dfrac{1}{k^2}+\dfrac{1}{\left(k+\dfrac{1}{k} \right)^2}}{k+\dfrac{1}{k+\dfrac{1}{k}}+\dfrac{1}{k}+\sqrt{k^2+4}}>0.
            \end{align*}
            From this we obtain that
            \begin{multline*}
                a_{n-i}+[0;a_{n-i+1},\dots,a_{n+2},\dots]+[0;a_{n-i-1},\dots,a_1]\\
                >\sqrt{k^2+4}+\dfrac{\dfrac{1}{k^2}+\dfrac{1}{\left(k+\dfrac{1}{k} \right)^2}}{k+\dfrac{1}{k+\dfrac{1}{k}}+\dfrac{1}{k}+\sqrt{k^2+4}}.
            \end{multline*}
            \item 
            Assume that $B>A$ and note that
            \begin{align*}
                &a_{n+1}+[0;a_{n+2},\dots]+[0;a_{n},\dots,a_{n-i-1},\dots,a_1]-\sqrt{k^2+4}\\
                &=k+[0;a_{n+2},\dots]+[0,\underbrace{k,k,\dots,k}_{i+1\text{ elements}},a_{n-i-1},\dots,a_1]-\sqrt{k^2+4}\\
                &>k+[0;a_{n+2},\dots]+[0,k,a_{n-i-1},\dots,a_1]-\sqrt{k^2+4}\\
                &=k+B+[0,k+A]-\sqrt{k^2+4}\\
                &>k+A+[0,k+A]-\sqrt{k^2+4}.
            \end{align*}
            Similarly, we also have that $f(A)=A+\dfrac{1}{k+A}$ is an increasing function for  $A>\dfrac{1}{k}$. Hence for all $A$ we have that
            $$f(A)>f\left( \dfrac{1}{k}\right)=\dfrac{1}{k+\dfrac{1}{k}}+\dfrac{1}{k}.$$
            From this we obtain in the same manner that
            \begin{align*}
                k+A+[0,k+A]-\sqrt{k^2+4}&>k+\dfrac{1}{k+\dfrac{1}{k}}+\dfrac{1}{k}-\sqrt{k^2+4}\\
                &=\dfrac{\dfrac{1}{k^2}+\dfrac{1}{\left(k+\dfrac{1}{k} \right)^2}}{k+\dfrac{1}{k+\dfrac{1}{k}}+\dfrac{1}{k}+\sqrt{k^2+4}}>0.
            \end{align*}
            This  implies that
            \begin{multline*}
                a_{n+1}+[0;a_{n+2},\dots]+[0;a_{n},\dots,a_{n-i-1},\dots,a_1]\\
                >\sqrt{k^2+4}+\dfrac{\dfrac{1}{k^2}+\dfrac{1}{\left(k+\dfrac{1}{k} \right)^2}}{k+\dfrac{1}{k+\dfrac{1}{k}}+\dfrac{1}{k}+\sqrt{k^2+4}}.
            \end{multline*}
        \end{enumerate}
    \end{enumerate}
\end{enumerate}
Therefore, the numbers have the worst approximation by rational
numbers in the form $\alpha=[a_0;a_1,\dots,a_m,\overline{k}]$, where $a_m\neq k$ and $m\geq1$. Hence for a large $n$ we have that
\begin{align*}
    \alpha-\dfrac{p_n}{q_n}&=\dfrac{(-1)^n}{q_n^2([\overline{k}]+[0;k,k,\dots,k,a_m,a_{m-1},\dots,a_1])}\\
    &=\dfrac{(-1)^n}{q_n^2(\sqrt{k^2+4}+\dfrac{k-\sqrt{k^2+4}}{2}+[0;k,k,\dots,k,a_m,a_{m-1},\dots,a_1])}.
\end{align*}
Set $r_m=[a_m;a_{m-1},\dots,a_1]$. This implies that
$$[0;k,k,\dots,k,a_m,a_{m-1},\dots,a_1]=[0;k,\dots,k,r_m]=\dfrac{p_{n-m}^*r_m+p_{n-m-1}^*}{q_{n-m}^*r_m+q_{n-m-1}^*}$$
where $\dfrac{p_j^*}{q_j^*}$ are convergents of $[0;\overline{k}]=\dfrac{\sqrt{k^2+4}-k}{2}$. From this and the fact that 
$$[0;\overline{k}]=\dfrac{p_{n-m}^*\dfrac{k+\sqrt{k^2+4}}{2}+p_{n-m-1}^*}{q_{n-m}^*\dfrac{k+\sqrt{k^2+4}}{2}+q_{n-m-1}^*}=\dfrac{p_{n-m}^*(k+\sqrt{k^2+4})+2p_{n-m-1}^*}{q_{n-m}^*(k+\sqrt{k^2+4})+2q_{n-m-1}^*}$$
we obtain that
\begin{align*}
    &\alpha-\dfrac{p_n}{q_n}=\dfrac{(-1)^n}{q_n^2\left(\sqrt{k^2+4}+[0;k,k,\dots,k,a_m,a_{m-1},\dots,a_1]-\dfrac{\sqrt{k^2+4}-k}{2}\right)}\\
    &=\dfrac{(-1)^n}{q_n^2\left(\sqrt{k^2+4}+\dfrac{p_{n-m}^*r_m+p_{n-m-1}^*}{q_{n-m}^*r_m+q_{n-m-1}^*}-\dfrac{p_{n-m}^*(k+\sqrt{k^2+4})+2p_{n-m-1}^*}{q_{n-m}^*(k+\sqrt{k^2+4})+2q_{n-m-1}^*}\right)}\\
    &=\dfrac{(-1)^n}{q_n^2\left(\sqrt{k^2+4}+\dfrac{1}{q_n^2}\dfrac{(k+\sqrt{k^2+4}-2r_m)(-1)^{n-m}q_n^2}{(q_{n-m}^*r_m+q_{n-m-1}^*)(q_{n-m}^*(k+\sqrt{k^2+4})+2q_{n-m-1}^*)}\right)},
\end{align*}
where \\
\begin{multline*}
    \dfrac{|k+\sqrt{k^2+4}-2r_m|q_n^2}{(q_{n-m}^*r_m+q_{n-m-1}^*)(q_{n-m}^*(k+\sqrt{k^2+4})+2q_{n-m-1}^*)} \\
    \geq \dfrac{|k+\sqrt{k^2+4}-2r_m|q_n}{q_{n-m}^*(k+\sqrt{k^2+4})+2q_{n-m-1}^*}.
\end{multline*}
From \eqref{worst} and the fact that $\sqrt{1+x}<1+\dfrac{1}{2}x$ for all $x\in(0,1)$ we obtain that for $\alpha=\dfrac{\sqrt{k^2+4}-k}{2}=[0;\overline{k}]$ or $\alpha=\dfrac{k+2-\sqrt{k^2+4}}{2}=[0;1,k-1,\overline{k}]$.   We have that
$$\left|\alpha-\dfrac{p_n}{q_n} \right|>\dfrac{1}{q_n^2\sqrt{k^2+4}+\dfrac{1}{\sqrt{k^2+4}}}.$$
So, to prove Theorem \ref{nathanson.t1}, it is enough to prove that if $\alpha\neq \dfrac{\sqrt{k^2+4}-k}{2}$ and $\alpha \neq \dfrac{k+2-\sqrt{k^2+4}}{2}$ plus an integer, then
$$\dfrac{|k+\sqrt{k^2+4}-2r_m|q_n}{q_{n-m}^*(k+\sqrt{k^2+4})+2q_{n-m-1}^*}>\dfrac{1}{\sqrt{k^2+4}},$$
for all sufficiently large $n$. Now the proof reduces to considering into several cases.
\begin{enumerate}
    \item 
    Assume that $a_m\geq k+1$. Then we have
    \begin{align*}
        &\dfrac{|k+\sqrt{k^2+4}-2r_m|q_n}{q_{n-m}^*(k+\sqrt{k^2+4})+2q_{n-m-1}^*}  \\
        &> \dfrac{(k+2-\sqrt{k^2+4})((k+1)q_{n-m}^*+q_{n-m-1}^*)}{q_{n-m}^*(k+\sqrt{k^2+4})+2q_{n-m-1}^*} >\dfrac{1}{\sqrt{k^2+4}}.
    \end{align*}
    
    \item 
    Let $a_m \in \{1,2,\dots,k-2 \}$, $(k\geq3)$ and $m \geq 1$. Then we obtain that
    \begin{align*}
        &\dfrac{|k+\sqrt{k^2+4}-2r_m|q_n}{q_{n-m}^*(k+\sqrt{k^2+4})+2q_{n-m-1}^*}  \\
        &> \dfrac{(\sqrt{k^2+4}-k+2)(q_{n-m}^*+q_{n-m-1}^*)}{q_{n-m}^*(k+\sqrt{k^2+4})+2q_{n-m-1}^*} >\dfrac{1}{\sqrt{k^2+4}}.
    \end{align*}
    
    \item 
    Suppose that $a_m=k-1, (k\geq2)$ and $m=1$. It implies that
    \begin{align*}
        &\dfrac{|k+\sqrt{k^2+4}-2r_m|q_n}{q_{n-m}^*(k+\sqrt{k^2+4})+2q_{n-m-1}^*}  \\
        &= \dfrac{(\sqrt{k^2+4}+2-k)((k-1)q_{n-m}^*+q_{n-m-1}^*)}{q_{n-m}^*(k+\sqrt{k^2+4})+2q_{n-m-1}^*} >\dfrac{1}{\sqrt{k^2+4}}.
    \end{align*}
    
    \item 
    Assume that $a_m=k-1, (k\geq2),a_{m-1}\geq2$ and $m\geq2$. Then it yields that
    \begin{align*}
        &\dfrac{|k+\sqrt{k^2+4}-2r_m|q_n}{q_{n-m}^*(k+\sqrt{k^2+4})+2q_{n-m-1}^*}  \\
        &\geq \dfrac{(\sqrt{k^2+4}+1-k)((2k-1)q_{n-m}^*+2q_{n-m-1}^*)}{q_{n-m}^*(k+\sqrt{k^2+4})+2q_{n-m-1}^*} >\dfrac{1}{\sqrt{k^2+4}}.
    \end{align*}
    
    \item 
    Let $a_m=k-1, (k\geq2),a_{m-1}=1$ and $m\geq3$. From this we obtain that
    \begin{align*}
        &\dfrac{|k+\sqrt{k^2+4}-2r_m|q_n}{q_{n-m}^*(k+\sqrt{k^2+4})+2q_{n-m-1}^*}  \\
        &\geq \dfrac{(\sqrt{k^2+4}-k)((2k-1)q_{n-m}^*+2q_{n-m-1}^*)}{q_{n-m}^*(k+\sqrt{k^2+4})+2q_{n-m-1}^*} >\dfrac{1}{\sqrt{k^2+4}}.
    \end{align*}
\end{enumerate}
The proof is complete.

\end{proof} 

\section{Data Availibility}

Data sharing is not applicable to this article as no new data were created or analysed in this study.

\section{Declaration}

The authors declare that they have no conflicts of interest. 

\section{Acknowledgement}
Tho Phuoc Nguyen is supported by grant SGS01/P\v{r}F/2024.

AMS Class: 11J82, 11A55.\\
Key words and phrases: continued fraction, approximation, Theorem of Nathanson.\\
Jaroslav Han\v{c}l, Tho Phuoc Nguyen, Department of Mathematics, Faculty of Sciences, University of Ostrava, 30.~dubna~22, 701~03 Ostrava~1, Czech Republic.\\
e-mail: jaroslav.hancl@seznam.cz, phuocthospt@gmail.com\\


\begin{thebibliography}{99}



\bibitem{bh} S. Bahnerov\' a, J. Han\v cl : Sharpening of the theorem of Vahlen and related theorems, J. Ramanujan. Math. Soc., vol. 36, no. 2, (2021), 109--121. 
\bibitem{borel1} \' E. Borel: Sur l'approximation des nombres par des nombres rationnels, C. R. Acad. Sci. Paris 136, (1903), 1054--1055. 
\bibitem{borel2} \' E. Borel: Contribution \`a l'analyse arithm\' etique du continu, J. Math. Pures 9, vol. 5, (1903), 329--375.
\bibitem{borwein} J. Borwein, P. Borwein, Pi and the AGM: A Study in Analytic Number Theory and Computational Complexity, John Wiley \& Sons, 
New York, 1987.
\bibitem{cohn} J. H. E. Cohn: Hurwitz's theorem. Proc. Amer. Math. Soc. 38 (1973), 436.
\bibitem{dirichlet} L. G. P. Dirichlet: Verallgemeinerung eines Satzes aus der Lehre von den Kettenbr\"uchen nebst einige Anwendungen auf die Theorie der Zahlen. S.-
B. Preuss. Akad. Wiss., (1842), 93–95.
\bibitem{feldman} N. I. Fel'dman, Yu. V. Nesterenko, Transcendental Numbers, Encyclopaedia of Mathematical Sciences, 
vol. 44: Number Theory IV, A. N. Parshin and I. R. Shafarevich, eds., Springer-Verlag, New York, 1998. 
\bibitem{hancl1} J. Han\v cl, Sharpening of theorems of Vahlen and  Hurwitz and approximation properties of the golden ratio, Arch. Math. (Basel) 105, no. 2, (2015), 129--137.  
\bibitem{hancl2} J. Han\v cl, Second basic theorem of Hurwitz, Lithuanian Mathematical Journal, vol. 56, no. 1, (2016), 72--76.  
\bibitem{hanclnair}  J. Han\v cl  and R. Nair: On a Theorem of Borel on Diophantine Approximation, Ramanujan J., (in press). 
\bibitem{hardy} G. H. Hardy,  E. M. Wright, An introduction to the theory of numbers. Sixth edition. Revised by D. R. Heath-Brown and 
J. H. Silverman. With a foreword by Andrew Wiles. Oxford University Press, Oxford, 2008. 
\bibitem{hensley} D. Hensley, Continued fractions, Word Scientific Publishing, (2006). 
\bibitem{hurwitz} A. Hurwitz, \"Uber die angen\" aherte Darstellung der Irrationalzahlen durch rationale Br\"uche, (German) Math. Ann. 39, 
no. 2, (1891), 279--284. 
\bibitem{jones} W. B. Jones, W. J. Thron, Continued fractions analytic theory and applications, Cambridge University 
Press, Encyclopedia of Mathematics and its applications 11, (1984). 
\bibitem{karpenkov} O. Karpenkov, \textit{Geometry of continued fractions.} Algorithms and Computation in Mathematics, 26. Springer, Heidelberg 2013. 
\bibitem{khinchin} A. Ya. Khinchin, Continued fractions, The University of Chicago Press, Chicago, (1964).  
\bibitem{nathanson} Nathanson M. B., Approximation by continued fractions, Proc Amer. Math. Soc. 45, (1974), 323--324. 
\bibitem{rosen} K. H. Rosen, Elementary number theory and its applications, Addison Wesley, fifth edition, (2005). 
\bibitem{schmidt} W. Schmidt, Diophantine approximation, Lecture Notes in Mathematics 785, Springer, Berlin, (1980). 
\bibitem{vahlen} K. Th. Vahlen, \" Uber N\" aherungswerte und Kettenbr\" uche, J. Reine Angew. Math. 115, (1895), 221--233.  
\bibitem{wall} H. S. Wall, Analytic theory of continued fractions, New York, Chelsea, (1948). 
\end{thebibliography}
\end{document}